\documentclass [12pt,leqno,a4paper]{article}
\usepackage{amsfonts}
\usepackage{amscd}
\catcode`\@=11
\def\resetthefootnote{\renewcommand{\thefootnote}{\@arabic\c@footnote} }
\def\@principiremex#1{\trivlist
  \item[\hskip \labelsep{\bfseries #1\ \thetheo}]\ignorespaces}
\def\opar@principiremex#1[#2]{\trivlist
  \item[\hskip \labelsep{\bfseries #1\ \thetheo\ (#2)}]\ignorespaces}

\newcommand{\newTHEOremrom}[2]{\newenvironment{#1}{\refstepcounter
  {theo}\@ifnextchar[{\opar@principiremex{#2}}{\@principiremex{#2}}
   }{\qedB\endtrivlist}}
\catcode`\@=12
\DeclareMathSymbol{\square}{\mathord}{AMSa}{"03}
\newcommand{\qed}{\nopagebreak\hspace*{\fill}{\vrule width6pt height6pt depth0pt}\par}
\newcommand{\qedB}{\nopagebreak\hspace*{\fill}$\square$\par}

\newif\ifAddress \Addressfalse

\newTHEOremrom{defi}{Definition}
\newTHEOremrom{rem}{Remark}
\newTHEOremrom{conj}{Conjecture}
\newTHEOremrom{ex}{Example}
\newTHEOremrom{nota}{Notation}
\newTHEOremrom{ques}{Question}

\newenvironment{trivlistparm}[2]{\trivlist \item[{#1 #2}]}{\endtrivlist}
\newenvironment{proclama}[1]{\trivlistparm{\bfseries}{#1}\itshape}{\endtrivlistparm}

\newenvironment{proof}{\trivlistparm{\itshape}{Proof.}}{\qed\endtrivlistparm}

\newcommand{\start}[2]{\begin{#1}\label{#2}}

\newcommand{\theoc}[1]{Theorem~\ref{#1}}
\newcommand{\propc}[1]{Proposition~\ref{#1}}
\newcommand{\coryc}[1]{Corollary~\ref{#1}}

\newcommand{\lemc}[1]{Lemma~\ref{#1}}

\newcommand{\remc}[1]{Remark~\ref{#1}}

\newcommand{\exc}[1]{Example~\ref{#1}}
\newcommand{\figc}[1]{Figure~\ref{#1}}

\newcommand{\refc}[1]{~\ref{#1}}
\newcommand{\refeq}[1]{\rlap(~\ref{#1})}

\newcounter{llistadepth}

\newenvironment{manlist}[1]{\addtocounter{llistadepth}{1}
      \edef\llistacontador{llista\romannumeral\the\value{llistadepth}}
      \list{({#1{\llistacontador}})}{\usecounter{\llistacontador}
      \def\makelabel##1{\hss\llap{##1}}
      \itemsep=2pt\parsep=0pt\topsep=3pt plus 1pt minus 1 pt}}{\endlist
      \addtocounter{llistadepth}{-1}}

\newenvironment{romlist}{\begin{manlist}{\roman}}{\end{manlist}}
\newenvironment{numlist}{\begin{manlist}{\arabic}}{\end{manlist}}

\def\dim{\mathrm{dim}}
\def\diag{\mathrm{diag}}

\def\wt{\widetilde}

\def\map#1#2#3{\mbox{${#1}\colon {#2} \longrightarrow {#3}$}}
\def\Smap#1#2{\mbox{${#1}\colon{#2} \longrightarrow {#2}$}}

\def\sett#1#2{\mbox{$\{{#1}\,\,:\,\,{#2}\}$}}

\newcommand{\Z}{\ensuremath{\mathbb{Z}}}
\newcommand{\Q}{\ensuremath{\mathbb{Q}}}
\newcommand{\C}{\ensuremath{\mathbb{C}}}
\newcommand{\R}{\ensuremath{\mathbb{R}}}
\newcommand{\Sc}{\ensuremath{\mathbb{S}}}
\newcommand{\SI}{\ensuremath{\mathbb{S}^1}}
\newcommand{\HH}{\ensuremath{\mathbb{H}}}

\def\hh{\mathop\mathcal{H}}

\def\s{\mathop\mathcal{S}}

\def\pow#1{\mbox{$(-1)^{#1}$}}
\def\powb#1{\mbox{$(-1)^{|{#1}|}$}}

\def\Map{\mathrm{Map}}
\def\mark{\mathrm{M \/}}
\def\erase{\mathrm{E \/}}

\def\lo{\Delta}

\def\flecha#1#2#3{\mbox{${#1}\stackrel{#2}{\longrightarrow} {#3}$}}

\def\flechi#1#2#3{\mbox{${#1}\stackrel{#2}{\mapsto} {#3}$}}

\def\hflecha#1#2{\mbox{$\stackrel{#1}{\longrightarrow} {#2}$}}

\newcommand{\ls}{\ensuremath{\mathbb{L}}}

\input epsf

\begin{document}

\title{String Topology}
\author{Moira Chas and Dennis Sullivan}
\maketitle

\begin{center}Abstract\end{center}
 Consider two families
of closed oriented curves in a manifold $M^d$. At
each point of intersection of a curve of one
family with a curve of the other family, form a
new closed curve by going around the first curve
and then going around the second. Typically, an
$i$-dimensional family and a $j$-dimensional
family will produce an $i+j-d+2$-dimensional
family. Our purpose is to describe a mathematical
structure behind such interactions.

\section{Introduction}
By the {\em string homology} of a manifold $M^d$
we mean the equivariant homology of the
continuous mapping space $\Map(\SI, M^d)$ with
the circle symmetry of rotating the domain. The
goal of this paper is to expose the following
theorem and its underpinnings:

\begin{proclama}{\theoc{gol}}
On the string homology of a smooth or
combinatorial oriented manifold $M^d$ there is a
natural graded Lie algebra structure of degree
$(2-d)$.
\end{proclama}

\begin{figure}[h]
\begin{center}
\epsfxsize=8cm \epsffile{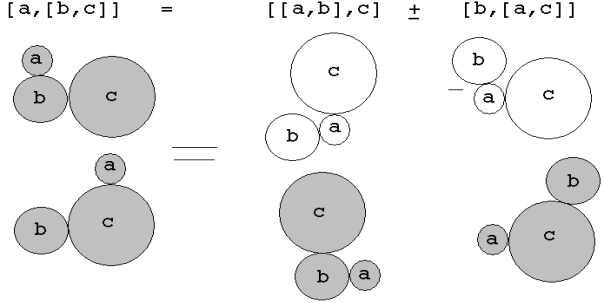}
\end{center}
\caption{Explanation of Jacobi identity
}\label{explanation}
\end{figure}

This structure called the {\em string bracket}
comes from the interaction of closed oriented
curves in $M^d$. At each point of intersection of
two such curves, one can form the oriented curve
obtained by going around one and then around the
other. Transversality will be used on the cells
computing string homology to define the string
bracket of the theorem. For an $i$ cell and a $j$
cell of strings there will typically be a
$(i+1)+(j+1)-d=(i+j)+(2-d)$ dimensional set of
intersection points. If we picture a bracket by
two circles touching, then a triple bracket means
a third circle touches the result along one arc
or the other. Thus, the Jacobi identity can be
viewed as in \figc{explanation}, where the two
terms on left appear on the right with two other
terms that are geometrically identical up to
sign. This Lie algebra is quite non-trivial for
surfaces of genus larger than one (see reference
to Goldman and Wolpert below).

In the process of analyzing the above argument we
found a more basic structure called the {\em loop
product} on the {\em loop homology}, the ordinary
homology of the free loop space $\Map(\SI, M^d)$.
A version of the loop product $\bullet$ of two
cells in $\Map(\SI,M)$ is defined for each point
$z$ on $\SI$ by first intersecting in $M^d$ the
two cells obtained by evaluating at this point
$z$ and then composing the loops at these
intersection points.

There are homotopies making the loop product on
loop homology into an associative graded
commutative algebra. As in Gerstenhaber's basic
paper \cite{Ger}, there is a preferred homotopy
for the graded commutativity (denoted $\ast$)
which leads by symmetrization to a second
operation $\{,\}$ on loop homology. This
operation, called {\em loop bracket} satisfies
the Jacobi identity where now it is convenient to
grade $\HH_\ast$ by subtracting $d$ from the
usual geometric grading so that $\bullet$ and
$\{,\}$ become operations of degree zero and one
respectively, $$ \flecha{\HH_i \otimes
\HH_j}{\bullet}{\HH_{i+j}}, $$ $$ \flecha{\HH_i
\otimes \HH_j}{\{,\}}{\HH_{i+j+1}}, $$ instead of
$-d$ and $(-d+1)$ respectively.

We arrive at the following result.

\begin{proclama}{\theoc{gersten}} The loop product $\bullet$ with the loop bracket
$\{,\}$ makes the loop homology $\HH_\star$ (the
ordinary homology of the free loop space) into a
Gerstenhaber algebra, namely:
\begin{numlist}
\item The loop product $\bullet$ defines a graded commutative, associative algebra.
\item $\{,\}$ is a Lie bracket of degree $1$, which means  that for each $a,b,c \in \HH_\ast$ \begin{romlist}
\item $\{a,b\}=-\pow{(|a|+1)(|b|+1)}\{b,a\}$
\item  $\{a,\{b,c\}\}=\{\{a,b\},c\}+\pow{(|a|+1)(|b|+1)}\{b,\{a,c\}\}$
\end{romlist}
\item $\{a,b\bullet c\}=\{a,b\}\bullet c+\pow{|b|(|a|-1)}b\bullet\{a,c\}$.
\end{numlist}
\end{proclama}

Now consider the circular symmetry of the loop
space $\Map(\SI, M)=\ls(M)$, and the associated
degree $+1$ operator $\Delta$ on homology
$$\flecha{i \mbox{ cycle in }\ls(M)}{}{i+1 \mbox{
cycle in } \ls(M)=i \mbox{ cycle}\times \SI}
\hflecha{\mbox{}{\scriptsize action}}{\ls(M)}. $$
It turns out that $\Delta$ is a second order
operator in the sense of commutative algebra,
i.e., the binary operation, the deviation of
$\Delta$ from being a derivation of $\bullet$, is
a derivation in each variable. One also has
easily that $\Delta \circ \Delta=0$. We come to

\begin{proclama}{\theoc{BV}} The loop product $\bullet$ and the operator $\Delta$ make  the loop homology (the ordinary homology of $\Map(\SI,M)$) into a Batalin Vilkovisky algebra, namely:
\begin{numlist}
\item $\bullet$ is a graded commutative associative algebra.
\item $\Delta \circ \Delta=0$.
\item $\pow{|a|}\Delta(a \bullet b)-\pow{|a|}\Delta a \bullet b-a \bullet \Delta b$ is a derivation of each variable.
\end{numlist}
\end{proclama}

The connection between these two results is that
conditions $(1)$, $(2)$, and $(3)$ of \theoc{BV}
imply formally that the binary operation defined
by the deviation satisfies graded Jacobi so that
a Batalin Vilkovisky algebra is a special type of
Gerstenhaber algebra. We prove \theoc{BV} by
constructing a chain homotopy between the
$\{x,y\}$ of \theoc{gersten} and
$\powb{x}\Delta(x\bullet y)-\powb{x}\Delta(x)
\bullet y-x\bullet \Delta(y)$. We note that in
the discussion of these two theorems there are
two independent proofs of Jacobi for the loop
bracket.

Finally we come again to the string bracket on
string homology, the equivariant homology of the
mapping space of $\SI$ into $M$. The approach
referred to above is to use intersection theory
on chains to define the string bracket and verify
the Jacobi identity directly on the level of
transversal triples of chains.

The second approach to the string bracket uses
the loop product $\bullet$ discussed above.
Consider the degree $+1$ operation {\it lift}
from equivariant chains to ordinary chains
corresponding to replacing an $i$-chain in the
base of an $\SI$ fibration by the $i+1$ chain
which is the preimage in the total space.
Consider also the operation {\it project} which
simply projects chains in the total space to the
base. Then we define the string bracket in terms
of the loop product by the formula $$
[x,y]=\mbox{project(lift $x \bullet$ lift $y$).}
$$ (In our calculations, Section \refc{3} we
denote {\em lift} by $\mark$, because lift means
mark a point in all possible ways on a closed
curve without a mark (a string). We denote {\em
project} by $\erase$ because project means erase
the marked point in a marked curve (loop) to get
an unmarked curve (string)).

The composition {\em (lift)$\bullet$ (project)}
induces the operator $\Delta$ above on loop
homology, and the Jacobi identity for $[,]$
follows by direct calculation from the properties
of the Gerstenhaber qua Batalin Vilkovisky
algebra structures above.

We arrive at \theoc{gol} stating that the string
bracket defines a graded Lie algebra structure on
string homology. The geometric degree is $(2-d)$
for the usual grading of string or equivariant
homology.

The same argument constructing the binary
operation $[,]=\bar{m}_2$, the string bracket,
constructs ternary, etc. operators $\bar{m}_3,
\bar{m}_4, \dots$.

Now extend each $\bar{m}_k$ to a coderivation
$m_k$ to $\Lambda \hh_\star$ (see
Section\refc{3}).

One knows the Jacobi identity for the string
bracket $\bar{m}_2$ is equivalent to the relation
$m_2 \circ m_2=0$ for the associated coderivation
$m_2$.

The Jacobi relation for $\bar{m}_2$ generalizes
to the entire collection $\{\bar{m}_2, \bar{m}_3,
\dots\}$ in the following way.

\begin{proclama}{\theoc{generalized}} The associated
coderivations $\{m_2, m_3, \dots\}$ of the free
commutative coalgebra $\Lambda \hh_\star$ on the
string homology $\hh_\star$ satisfy:
\begin{romlist}
\item $m_k \circ m_k=0$, for $k=2,3,4\dots$.
\item $m_k \circ m_r + m_r \circ m_k=0$ for $k,r=2,3,4,\dots$.
\end{romlist}
\end{proclama}

These coderivations $m_k$ combine in various ways
to define coderivations of square zero on
$\Lambda\hh_\star$. Such a differential is one
definition of a $\mathrm{Lie}_\infty$ or strong
homotopy Lie algebra structure on $\hh_\star$.
Thus we have

\begin{proclama}{\coryc{uncountable}}
There exists an uncountable family
$\{\delta_\Lambda\}$ of $\mathrm{Lie}_\infty$
structures on the string homology. Namely, for
each $\Lambda \subset \{2,3,\dots\}$,
$$\Smap{\delta_\Lambda}{\Lambda \hh_\star}\mbox{
defined as } \delta_\Lambda=\sum_{\lambda \in
\Lambda}m_\lambda $$ is a coderivation which
satisfies $\delta_\Lambda \circ
\delta_\Lambda=0$.
\end{proclama}

If we examine the string bracket when $d=2$ we
find a Lie bracket structure on the vector space
of components of the space of closed curves in a
surface. For surfaces of genus larger than zero
this is the non-trivial bracket discovered in the
80's by Wolpert \cite {Wol} and Goldman
\cite{Gol}. That discovery was strongly related
to the symplectic structure on Techmuller space
\cite{Wol} and the symplectic structure of other
spaces of flat connections over a surface
\cite{Gol}. (For more discussion see
Section\refc{examples}). When the genus is zero,
i.e., $M$ is a two sphere, the loop product
becomes non-trivial in the higher dimensional
algebraic topology of the free loop space of
$\Sc^2$.

The $\Sc^2$-calculation is part of a general
structure based on the diagram relating the
intersection product on ordinary homology
$(H_\star(M),\wedge)$ with the loop product on
loop homology $(\HH_\star, \bullet)$ and the
Pontryagin product on the based loop space
homology $(H_\star(\Omega), \times)$,

\begin{equation}
\flecha{(H_\star(M),\wedge)}{\mbox{{\scriptsize
constant loops}}}{(\HH_\star,
\bullet)}\hflecha{\mbox{{\scriptsize intersection
with a fiber}}}{(H_\star(\Omega),\times)}
\label{diagram}
\end{equation}

Both maps are ring homomorphisms. The first is an
injection onto a direct summand showing {\it the
loop product is an extension of the classical
intersection product}. The image of the second
map is a {\it graded commutative} subalgebra. For
$\Sc^2$, $(H_\star(\Omega),\times)$ is the tensor
algebra on one generator $\eta$ in degree one and
the image is the subalgebra generated by
$\eta^2$. This shows the loop product is
non-trivial for $\Sc^2$ (Section \refc{two
sphere}). For more complete calculations we can
augment the diagram \refeq{diagram} with a
relation between the usual cap product operation,
$\cap$, and the loop product. This is described
by

\begin{proclama}{\theoc{compatible}} For each $x, y$ homology classes and
compatible pair of classes $(A,a)$ $$ \mbox{loop
product }( A \cap(x \otimes y))=a \cap(x \bullet
y). $$
\end{proclama}
(see Section \refc{cup} for the definitions and
the proof, and Section \refc{two sphere} for
relevant algorithms.)

The two approaches here to the string bracket,
direct geometry and via the loop product, reminds
one of Witten's paper \cite{Witt}. There it was
pointed out that closed string interactions
looked at directly as in the \figc{explanation}
are non-associative. To get around this, a
marking point was introduced in \cite{Witt} to
facilitate the definition of an associative
multiplication of (open) strings. Thus also our
string bracket is non-associative but satisfies
Jacobi and it arises from an associative product
of loops (marked strings).

There is also a dictionary relating our
constructions with those in algebra begun by
Gerstenhaber \cite{Ger}.

In a sequel we will discuss a rich world of
general string operations in the chains of the
loop space. We find a structure like a big part
of a two dimensional field theory associated to
each manifold $M^d$. In particular we investigate
the $\mathrm{Lie}_\infty$ structure described in
\theoc{generalized} and \coryc{uncountable} as
well as co-versions. We also hope to follow
Stasheff's specific suggestion to relate our
structure to the work of the physicist Zwiebach
\cite{WZ} and \cite{Z}.

\begin{center}{\bf Table of contents}
\end{center}
\begin{numlist}
\item Introduction
\item Definition of the loop product
$\bullet$ at the chain level and passage of
$\bullet$ to loop homology.
\item $\bullet$ is homotopy associative and homotopy commutative where  $\ast$ is the homotopy of commutativity.
\item Symmetrization of $\ast$ gives a loop bracket $\{,\}$ which passes to loop homology and we have a Gerstenhaber algebra (i.e., the bracket is a biderivation of $\bullet$ and satisfies Jacobi).
\item Definition of  $\Delta$ and the argument that the deviation of
$\Delta$ from being a derivation of $\bullet$ is
the loop bracket $\{,\}$.
\item Definition of  $\mark$ and $\erase$,  string homology and   the construction of a string bracket  $[,]$ which satisfies
Jacobi. Generalized $n$-brackets and generalized
Jacobi.
\item The string bracket for surfaces, the work of Wolpert and Goldman, and the case of $\Sc^2$.
\item Cap product and loop product.
\item Appendix 1: $Sc^2$ and other simply
connected manifolds.
\item Appendix 2: $M^3$ and $K(\pi,1)$
manifolds.
 \end{numlist}

\section{The loop product $\bullet$}\label{1}

We think of the circle $\SI$ as $\R/\Z$, so each
point can specified by some $x \in [0,1)$. Unless
otherwise stated, $M$ is an orientable manifold
of dimension $d$. A {\em loop in $M$} is a
continuous map from $\SI$ to $M$. Observe that a
loop has a marked point: the image of $0$.
$\Map(\SI, M)$, the space of all loops, will be
denoted by $\ls(M)$, or by $\ls$ since $M$ is
fixed throughout the discussion. By an $i$-chain
we mean a linear combination of oriented
$i$-dimensional families of loops in $M$. The
parameter spaces of the families are taken from
any standard list of cells closed under face
operators.

In the algebraic topology of chains on a space
there are two well known multiplications which we
will combine to obtain a new structure. The first
of these multiplications is the (transversal)
intersection of chains in a $d$-manifold: an
$i$-chain intersected with a $j$-chain gives an
$i+j-d$-chain.

The second is the product of an $i$-chain of
loops with a $j$-chain of loops, all of whose
marked points are equal to some $p \in M$.
Multiplying these chains yields an $i+j$
cartesian product chain of composed loops with
marked point at $p$ in $M$.

\begin{figure}[h]
\begin{center}
\epsfxsize=10cm \epsffile{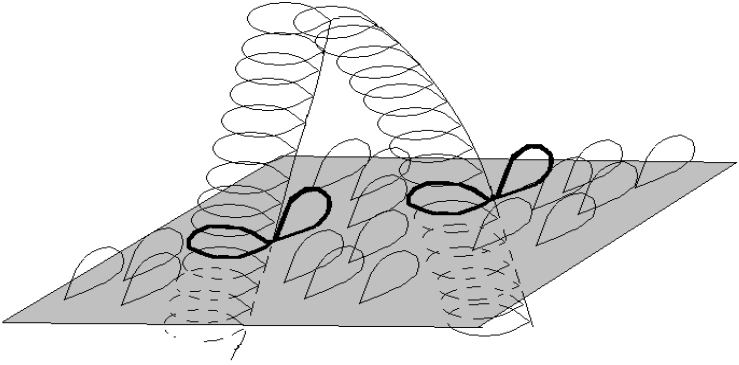}
\end{center}
\caption{ For $M^3$, the loop product of a
$1$-cell with a $2$-cell is a
$0$-chain}\label{products}
\end{figure}

Our new loop product $\bullet$ is transversally
defined at the chain level as follows (see
\remc{dos}): given $x$, an $i$-chain of loops in
$M$, and $y$, a $j$-chain of loops in $M$, one
first intersects the $i$-chain (of $M$) of marked
points of $x$ with the $j$-chain (of $M$) of
marked points of $y$, to obtain an $i+j-d$ chain
$c$ (of $M$) along which the marked point of $x$
coincides with the marked point of $y$. Now we
define the chain $x \bullet y$, by putting at
each point of $c$ the composed loop that first
goes around the loop of $x$ and then, around the
loop of $y$ (see \figc{products}).

In order to give the precise definition of the
loop product we need the following remark.
\start{rem}{uno} We will use the following {\bf
orientation convention} for our constructions. We
will have a map of a domain cell, usually a
product $K_1\times K_2\times \dots K_n$ into
$M\times M \times \dots M$ ($l$ factors) which is
assumed to be transversal to a diagonal. We
orient the product cell with the product of
individual orientations. We orient the normal of
the pull back by pulling back the orientation of
the normal of the diagonal induced by some
orientation of $M$. We then take the induced
orientation on the pullback so that $$
\mbox{(orientation on the pullback)(orientation
on the normal)}=$$$$\mbox{(product orientation of
the cell domain)} $$
\end{rem}

Denote the $(i-d)$ chains of $\ls$ by $\ls_i$ and
the direct sum of all these by $\ls_\ast$. Also,
if $\map{x}{K_x}{\ls}$ is a cell, we denote by
$K_x$ its underlying set.

Now, for any pair of cells, $\map{x}{K_x}{\ls}$,
$\map{y}{K_y}{\ls}$ we define the set $K_{x
\bullet y}$ as the transversal preimage of the
diagonal of $M \times M$ under the map
$$\flecha{K_x \times K_y}{}{M \times M},$$
$$(k_x,k_y) \mapsto (x(k_x)(0),y(k_y)(0)) $$ Now,
define $\map{x\bullet y}{K_{x \bullet y}}{\ls}$
as
\begin{eqnarray*}
(x\bullet y)(k_x,k_y)(\gamma) &=
&\left\{\begin{array}{ll} x(k_x)(2\gamma) &
\mbox{if $\gamma \in [0,\frac{1}{2}]$,}\\
y(k_y)(2\gamma) & \mbox{if $\gamma \in
[\frac{1}{2},1)$.}
\end{array} \right.
\end{eqnarray*}
(we keep the notation $K_{x\bullet y}$ because
$K_{x \bullet y}$ is a manifold which can be
divided into cells) . Orient $K_{x\bullet y}$ by
\remc{uno}. (Observe that with our new grading,
if $x \in \ls_i$ and $y \in \ls_j$ then $x
\bullet y \in \ls_{i+j}$.

\start{rem}{dos} We will adopt a point of view
which is used in classical intersection theory of
chains in a manifold. We say that a chain
operation is {\em transversally defined} if it is
defined for appropriately transversal cells.

We say that an identity between chain operations
{\em holds transversally} if it holds on any
finite subset of the chains where all
constituents are appropriately transversal. Thus,
the classical intersection is defined
transversally at the chain level, is
transversally associative and transversally
graded commutative.
\end{rem}

\start{lem}{passhom} If $x, y \in \ls_\ast$ is a
transversal pair then $$\partial(x \bullet
y)=\partial x \bullet y +\pow{|x|}x \bullet
\partial y.$$
\end{lem}
\begin{proof}  By definition, the underlying chain of $x \bullet y$ is the oriented intersection chain in $M$ of the marked points of $x$ and of $y$ respectively. Once the orientation of $M$ is fixed, the orientations of intersections behave for calculations like orientations of normal directions.

Now, for each $x \in \ls_\ast$, we denote by
$\wt{x}$ the chain of $M$ of marked points of
$x$. Thus, for transversal intersection one gets
the familiar $$
\partial(\wt{x}\cap \wt{y})=(\partial \wt{x} \cap y)+\powb{x}(\wt{x} \cap \partial \wt{y})
$$ where $|x|=\dim(x)-d$.

Since $K_{x \bullet y}$ is the underlying set of
$\wt{x}\cap \wt{y}$ we get $$
\partial K_{x \bullet y}=K_{\partial x \bullet y}+\powb{x}K_{x \bullet \partial y}.
$$

The chain $\partial(x \bullet y)$ is the
restriction of $x \bullet y$ to $\partial
K_{x\bullet y}.$ Therefore, the above formula
yields the same formula for $\bullet$, $$
\partial(x\bullet y)=\partial x \bullet y +\pow{|x|} x \bullet \partial y.
$$
\end{proof}

Let $\HH_i$ denote the $i$-th homology group of
the loop space $\ls$ with the degree shifted down
by $d$ and set $\HH_\star=\HH_{-d} \oplus
\HH_{-d+1}\oplus \dots \HH_0 \oplus \HH_1\dots $.
We refer to $\HH_\star$ as the {\em loop
homology}.

\start{cory}{tres} The loop product $\bullet$ on
chains passes to loop homology and defines a
product $$ \flecha{\HH_i \otimes
\HH_j}{\bullet}{\HH_{i+j}}. $$

\end{cory}

\start{rem}{cuatro} The loop product is defined
in all dimensions and may be non-zero way above
the dimension of the manifold (e.g., $M=\Sc^2$,
\exc{gold}, or $\Sc^3$, \coryc{sphere}).
\end{rem}

\section{Associativity of $\bullet$, the $\ast$ operator, and homotopy commutativity of $\bullet$}\label{2}

Let us first discuss associativity. For the
classical intersection product, we have that if
three cycles are pairwise transversal, then the
intersection product is literally associative at
the chain level.

The classical based loop composition is
associative up to homotopy (see Stasheff
\cite{Sta1} and \cite{Sta2} for a complete
discussion of this point).

Thus both of these classical chain products yield
associative multiplications in homology, and thus
we will have the same associativity for
$\bullet$, the loop product, combining
intersection and based loop composition.

\start{prop}{associativity} The loop product in
loop homology $\flecha{\HH_i \otimes
\HH_j}{\bullet}{\HH_{i+j}}$ is associative.
\end{prop}
\begin{proof} Assume that the three homology classes in question $x, y, z$
are represented by cycles that are pairwise
transversal. The intersection locus of $(x\bullet
y)\bullet z$ and $x \bullet (y \bullet z)$ are
literally equal with identical coorientations.
The loop product is associative up to homotopy
using the same considerations as in the based
loop product, now parameterized by the points of
the intersection $x \cap y \cap z$.
\end{proof}

\begin{figure}[h]
\begin{center}
\epsfxsize=5cm \epsffile{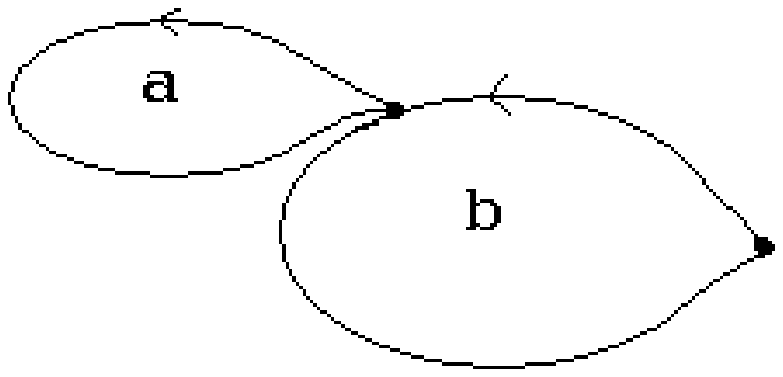}\end{center}
\caption{The $\ast$ product }\label{333}
\end{figure}

Let us turn to the question of commutativity. The
intersection product of two transversal cycles is
literally graded commutative. However, the based
loop product (Pontryagin product) is often
non-commutative even in homology. We will see
that the $\bullet$ product, combining these two
products is homotopy commutative at the chain
level. One chain homotopy is given by a new
binary operation $x \ast y$ defined for
appropriately transversal pairs $x, y$ in the
following way. Consider the chain $c$ (of $M$)
where the marked point of $x$ transversally
intersects one of the images of the loops of $y$.
Then at each point of $c$ put the following loop:
first go around the loop of $y$ until the
intersection point with the marked point of $x$.
Then go around $x$ and finally, go around the
rest of $y$ (see \figc{333}).

The precise definition is the following: let
$\map{x}{K_x}{\ls}$, $\map{y}{K_y}{M}$ be two
cells in $\ls_\star$. Let $K_{x \ast y}$ be the
preimage of the diagonal of $M \times M$ under
the map $$\flecha{K_x \times [0,1] \times
K_y}{}{M \times M},$$ $$(k_x,s,k_y) \mapsto
(x(k_x)(0),y(k_y)(p(s))),$$ where
$\map{p}{[0,1]}{\SI}$ is the usual projection.
Then $\map{x\ast y}{K_{x \ast y}}{\ls}$ is
defined as

\begin{eqnarray*}
(x\ast y)(k_x,s,k_y)(\gamma) &=
&\left\{\begin{array}{ll} y(k_y)(2\gamma) &
\mbox{if $\gamma \in [0,\frac{s}{2})$.}\\
x(k_x)(2\gamma-s) & \mbox{if $\gamma \in
[\frac{s}{2},\frac{s+1}{2})$,}\\ y(k_y)(2\gamma)
& \mbox{if $\gamma \in [\frac{s+1}{2},1)$.}
\end{array} \right.
\end{eqnarray*}

\start{lem}{lemma2} If $x, y $ $\in \ls_\star$
are appropriately transversal, then $$\partial(x
\ast y)=\partial x \ast y+\pow{|x|+1}x \ast
\partial y+\pow{|x|}(x\bullet y
-\pow{|x||y|}y\bullet x). $$
\end{lem}
\begin{proof} Let $\wt{x}$ denote the chain of $M$ of the marked points of $x$ and let ${\wt{y}}$ denote the chain of $M$ with parameter space $[0,1] \times K_y$, given by $\wt{y}(s,k_y)=y(k_y)(p(s))$, where $\map{p}{[0,1]}{\SI}$ is the usual projection.

Suppose that $x$ and $y$ are appropriately
transversal. Since the underlying sets of $x \ast
y$ and $\wt{x} \cap \wt{y}$ coincide, $\partial(x
\ast y)$ is parametrized by the underlying set of
$\partial(\wt{x} \cap \wt{y})$. Hence, for
transversal intersection, one obtains $$
\partial(\wt{x}\cap \wt{y})=(\partial \wt{x} \cap y)+\powb{x}(\wt{x} \cap \partial \wt{y}).\label{hi}
$$ Since $\partial([0,1]\times K_y)=\partial[0,1]
\times K_y-[0,1]\times \partial K_y$, the
restriction of $x \ast y$ to the underlying set
of $\powb{x}(\wt{x} \cap \partial \wt{y})$ is
$\pow{|x|}(x\bullet y -\pow{|x||y|}y\bullet
x)-\powb{x}x \ast
\partial y.$ One the other hand, the restriction of $x \ast y$ to
the underlying set of $\partial \wt{x} \cap y$ is
${\partial x \ast y}$, which completes the proof.
\end{proof}

In particular, if $x$ and $y$ are cycles,
\lemc{lemma2} implies that $\partial(x \ast
y)=\pm(x\bullet y -\pow{|x||y|}y\bullet x).$
Therefore, we have

\start{theo}{lemma1} $(\HH_\star, \bullet)$ is an
associative, (graded) commutative algebra.
\end{theo}

Let us compare the loop product on loop homology
$(\HH_\star, \bullet)$ with the usual homology of
the manifold with intersection product
$(H_\star(M), \wedge)$ and the homology of the
based loop space with the based loop product or
Pontryagin product, $(H_\star(\Omega), \times)$.
We have two maps, $$
\flecha{H_\star(M)}{\varepsilon}{\HH_\star}\hflecha{\cap}{H_\star(\Omega)}
$$ where $\varepsilon$ is the inclusion of
constant (or even $\varepsilon$ small) loops into
all loops and $\cap$ is the transversal
intersection with one fiber of the projection
loop space $\hflecha{\mbox{evaluation}}{M}$. If
we use the usual grading on $H(\Omega)$, our
shifted grading on $\HH_\star$, the homology of
the entire loop space, and the analogous shifted
grading on $H_\star(M)$, then these products and
the two maps have degree zero.

\start{prop}{maps} $\flecha{(H_\star(M),
\wedge)}{\varepsilon}{(\HH_\star,
\bullet)}\hflecha{\cap}{(H_\star(\Omega),\times)}$
preserve products.\end{prop}

\begin{proof}These follow directly from the definitions.
\end{proof}
\start{rem}{injection} $\varepsilon$ is an
injection onto a direct summand. For any Lie
group manifold, $\cap$ is a surjection.
\end{rem}

\start{cory}{sphere} For $M=\Sc^3$, the loop
product is non-zero in infinitely many degrees.
\end{cory}
\begin{proof}By \remc{injection}, $\flecha{\HH_\star}{\cap}{H_\star(\Omega)}$ is a surjection. With $\Q$ coefficients, the homology of $\Omega(\Sc^3)$ is the homology of $\C P^\infty$
and the Pontryagin product gives a polynomial
algebra on this generator in degree $2$.
\end{proof}

\section{The loop bracket $\{,\}$ in loop homology}\label{33}

For easier thinking, let us now derive formulae
by calculating in the chain complex of all
homomorphisms $\flecha{\ls_\ast \otimes
\ls_\ast}{\varphi}{\ls_\ast}$ of degree $-1,
0,1,\dots$, with the usual $\partial$ defined by
the Leibniz rule with signs, $$
\partial(\varphi(a\otimes b))=(\partial \varphi)(a\otimes b)+\powb{\varphi}\varphi(\partial(a\otimes b)),
$$ $\partial \varphi$, which is defined by this
relation, is also denoted $[\partial,\varphi]$,
the graded commutator of $\partial$ and
$\varphi$.

If $[\partial,\varphi]=0$, $\varphi$ is usually
called a {\em chain map}. As an example, if
$\varphi(x \otimes y)=\powb{x}x\ast y$ of the
previous section, $[\partial,\varphi]$ evaluated
on $(x\otimes y)$ is $$\powb{x}\partial(x \ast
y)+(\pow{|x|+1}\partial x \ast y+x \ast \partial
y)=$$ $$ \powb{x}(\partial(x\ast y)-\partial x
\ast y -\pow{|x|+1}x\ast
\partial y).$$ By \lemc{lemma2} this is equal to $$x\bullet y
-\pow{|x||y|}y\bullet x.$$

In other words, in the chain complex of
homomorphisms $\flecha{\ls_\ast \otimes
\ls_\ast}{\varphi}{\ls_\ast}$, the transversally
defined $1$-chain $\powb{x}x\ast y$ has boundary
the transversally defined zero chain $(x\bullet y
-\pow{|x||y|}y\bullet x).$ Symbolically, $$
\partial \ast'=\bullet-\bullet \tau, $$ where $x
\ast' y=\powb{x}x \ast y$ and $\tau(x\otimes
y)=\pow{|x||y|}y \otimes x$.

Thus we consider $\ast'+\ast'\tau $ and calculate
it is a $1$-cycle $$[\partial,\ast'+\ast'\tau
]=[\partial,\ast']+[\partial,\ast'\tau]=$$ $$
[\partial,\ast']+[\partial,\ast']\tau=(\bullet-\bullet
\tau)+(\bullet \tau-\bullet)=0 $$ using
$[\partial, \tau]=0$ and $\tau^2=0$.

\begin{defi} The loop bracket $\{x,y\}$ is defined transversally on $\ls_\ast$ by the formula
$$ \{x,y\}=x\ast y-\pow{(|x|+1)(|y|+1)}y\ast x.
$$
\end{defi}

\start{lem}{lemma3} For $x,y,z \in \ls_\ast$ the
associator of $\ast$ is symmetric in the first
two variables, $$ x \ast (y \ast z)-(x \ast
y)\ast z=\pow{(|x|+1)(|y|+1)}(y\ast (x\ast
z)-(y\ast x)\ast z). $$
\end{lem}
\begin{proof} Observe that the parameter spaces of $ K_{x \ast (y \ast z)-(x \ast y)\ast z}$ and $K_{y\ast (x\ast z)-(y\ast x)\ast z}$ consist in those values where the basepoint of
the loops of $x$ coincides with the one of the
points of the loop of $z$ and the basepoint of
loops of $y$ coincides with another. (see
\figc{44})

\begin{figure}[h]
\begin{center}
\epsfxsize=11cm \epsffile{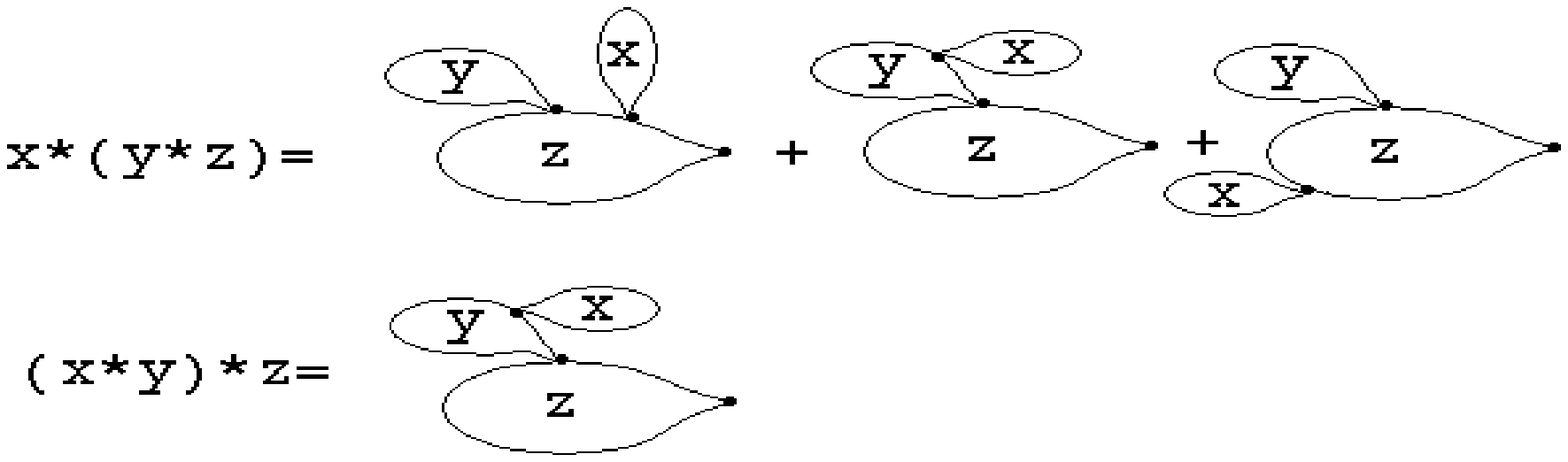 }
\end{center}
 \caption{ $K_{x\ast(y\ast z)}$ and
$K_{(x\ast y)\ast z}$}\label{44}
\end{figure}

Symbolically, $$ K_{x \ast (y \ast z)-(x \ast
y)\ast z}=\{(k_x,s,k_y,t,k_z) \in K_x \times
[0,1] \times K_y \times [0,1] \times
K_z:$$$$x(x_x)(0)=z(k_z)(s),
y(k_y)(0)=z(k_z)(t)\} $$ $$K_{y\ast (x\ast
z)-(y\ast x)\ast z}=\{(k_y,t,k_x,s,k_z) \in K_x
\times [0,1] \times K_y \times [0,1] \times
K_z:$$$$x(x_x)(0)=z(k_z)(s),
y(k_y)(0)=z(k_z)(t)=\}. $$ Clearly, there is a
bijection between the underlying sets of the
chains in question, $K_{x \ast (y \ast z)-(x \ast
y)\ast z}$ and $K_{y\ast (x\ast z)-(y\ast x)\ast
z}$, which is orientation preserving if and only
if $\pow{(|x|+1)(|y|+1)}=1$.
\end{proof}

\start{prop}{Jacobi} $(\ls_\star,\{,\})$ is a
graded Lie algebra (transversally) with all the
degrees shifted by $1$. In other words, for each
$x, y, z \in \ls_\star$ mutually transversal,
\begin{numlist}
\item $\{x,\{y,z\}\}=\{x,\{y,z\}\}+\pow{(|x|+1)(|y|+1)}\{y,\{x,z\}\}.$
\item $\{x,y\}=-\pow{(|x|+1)(|y|+1)}\{y,x\}$
\end{numlist}
\end{prop}
\begin{proof}Let us prove $(1)$. By \lemc{lemma3},
$$\{\{x,y\},z\}+\pow{(|x|+1)(|y|+1)}\{y,\{x,z\}\}-\{x,\{y,z\}\}=$$
$$(x \ast y)\ast z - \pow{(|x|+1)(|y|+1)}(y \ast
x)\ast z -\pow{(|x|+|y|)(|z|+1)}z \ast (x \ast y)
+\pow{|x||z|+|y||z|+|x||y|+1)}z\ast (y \ast x)
$$$$ +\pow{(|x|+1)(|y|+1)}y\ast(x\ast z)
-\pow{(|x|+1)(|y|+|z|)}y\ast(z\ast x)$$$$
-\pow{(|y|+1)(|z|+1)}(x\ast z)\ast y
+\pow{(|x|+|y|)(|z|+1)}(z\ast x)\ast y $$$$-x\ast
(y\ast z) +\pow{(|y|+1)(|z|+1)}x\ast (z \ast y)
+\pow{(|x|+1)(|y|+|z|)}(y\ast z)\ast x$$$$
-\pow{|x||y|+|y||z|+|x||z|+1}(z\ast y)\ast x =0$$
\end{proof}

\start{rem}{associator} An identical Proposition
and proof can be found in \cite{Ger} in a purely
algebraic context.
\end{rem}

\start{cory}{jac} The loop homology with the loop
bracket,$(\HH,\{,\})$ is a graded Lie algebra of
degree $+1$.
\end{cory}

Now we discuss the compatibility of loop bracket
and loop product.

By \lemc{lemma2}, the map (transversally defined)
$$\flecha{\ls_\star \otimes \ls_\star\otimes
\ls_\star}{}{\ls_\star},$$ $$\flecha{x_1 \otimes
x_2 \otimes y}{}{(x_1 \bullet x_2)\ast y -x_1
\bullet (x_2 \ast y)-\pow{|x_1|(|y|+1)}(x_1 \ast
y) \bullet x_2}$$ is a chain map, i.e., the
commutator with $\partial$ is zero transversally.

\start{lem}{star} Let $x, x_1, x_2, y, y_1, y_2
\in \ls_\ast$ be appropriately transversal, then
\begin{numlist}
\item $x\ast (y_1 \bullet y_2)=(x \ast y_1)\bullet y_2+\pow{|y_1|(|x|+1)}y_1 \bullet (x\ast y_2)$
\item $(x_1 \bullet x_2)\ast y -x_1 \bullet (x_2 \ast y)-\pow{|x_1|(|y|+1)}(x_1 \ast y) \bullet x_2 \simeq 0$, where $\simeq$ means chain homotopy.
\end{numlist}
\end{lem}
\begin{proof}
The proof of $(1)$ is easier because the equation
holds transversally at the chain level: The set
of parameters where the marked point of $x$
coincides with one of the images of $y_1 \bullet
y_2$ is the union of the set of parameters where
the marked point of $x$ coincides with one of the
image of $y_1$ union the set of parameters where
the marked point of $x$ coincides with one of the
image of $y_2$. (see \figc{55})

\begin{figure}[h]
\begin{center}
\epsfxsize=12cm \epsffile{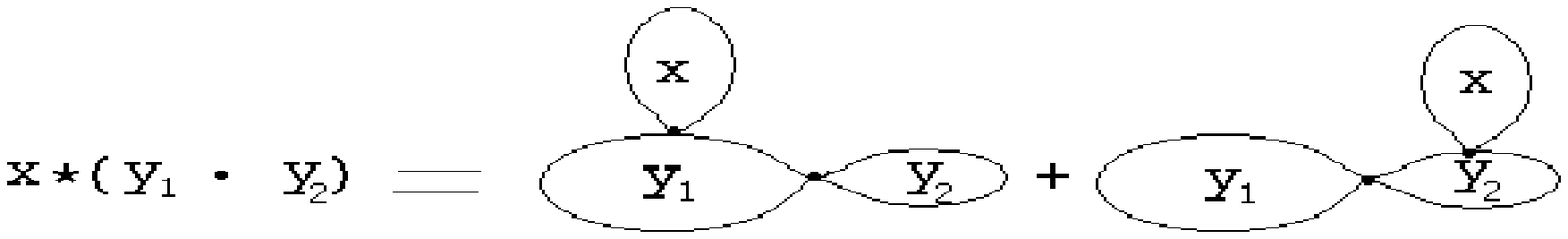}
\end{center}
\caption{Proof of \lemc{star}(1)}\label{55}
\end{figure}

Let us prove $(2)$. The idea (see \figc{66}) is
that the $2$ dimensional chain
$$\map{\varphi}{\ls_\ast \otimes \ls_\ast \otimes
\ls_\ast}{\ls_\ast},$$ $$\flecha{x_1 \otimes x_2
\otimes y}{\varphi}{\varphi_{x_1,x_2,y}}$$ where
$x_1$ and $x_2$ attach to $y$ at pairs of
arbitrary points in such a way relative to the
cyclic order that $x_1$ is between the marked
point and $x_2$, provides a chain homotopy
between the two sides. More precisely, for each
pair of points $$(s,t) \in T=\sett{(s,t) \in
[0,1]\times [0,1]}{s +t \le 1},$$ we will define
a chain such that $x_1$ is attached to $y$ at $t$
and $x_2$ is attached to $y$ at $1-s$.

Let $$K=\{(k_{x_1},k_{x_2},s,k_y,t)\in
K_{x_1}\times K_{x_2}\times [0,1]\times K_y
\times [0,1]:$$ $${x_1}(k_{x_1})(0)=y(k_y)(t),
{x_2}(k_{x_2})(0)=y(k_y)(s), s,t\in [0,1];s+t\le
1\}$$

We define a map
$\map{\varphi_{x_1,x_2,y}}{K}{\ls}$.
\begin{eqnarray*}
\varphi_{x_1,x_2,y}(k_{x_1},k_{x_2},s,k_y,t)(\gamma)
&= &\left\{\begin{array}{ll} y(3\gamma) &
\mbox{if $\gamma \in [0,\frac{t}{3}]$}\\
{x_1}(3\gamma-t) & \mbox{if $\gamma \in
[\frac{t}{3}, \frac{t+1}{3}]$.}\\ y(3\gamma) &
\mbox{if $\gamma \in
[\frac{t+1}{3},\frac{1}{3}]$}\\
{x_2}(3\gamma+s-1) & \mbox{if $\gamma \in
[\frac{1-s}{3},\frac{2-s}{3}]$.}\\ y(3\gamma) &
\mbox{if $\gamma \in [\frac{2-s}{3},1]$}.
\end{array}
\right.
\end{eqnarray*}
for each $(k_{x_1},k_{x_2},s,k_y,t) \in K$,
$\gamma \in \SI$.

\begin{figure}[h]
\begin{center}
\epsfxsize=10cm \epsffile{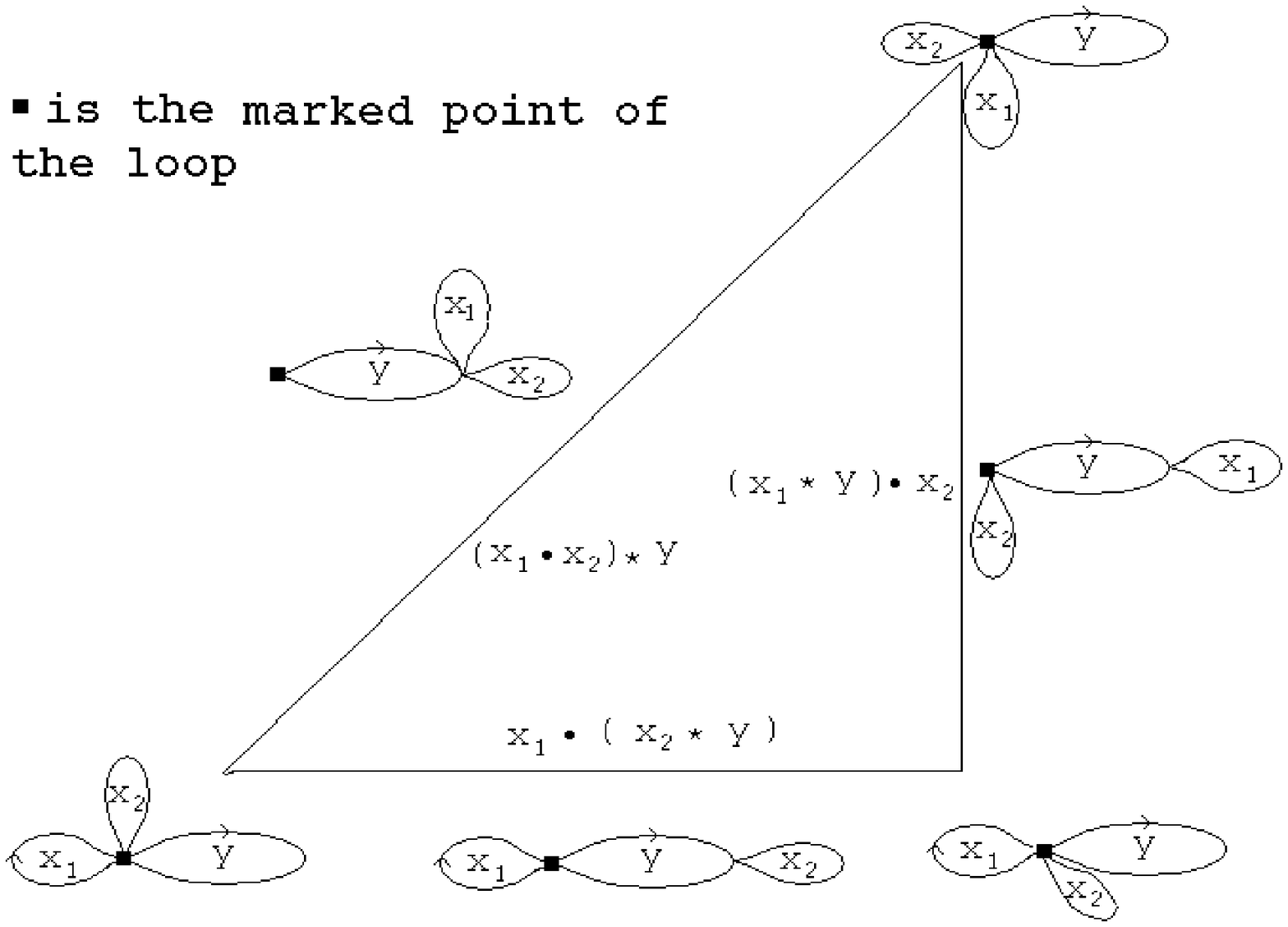}
\end{center}
\caption{Proof of \lemc{star}(2)}\label{66}
\end{figure}

Since the restriction of $\varphi$ to
\begin{romlist}
\item $K \cap \sett{(k_{x_1},k_{x_2},s,k_y,t) \in K_{x_1}\times K_{x_2}\times [0,1]\times K_y \times [0,1]}{t=0}$ is $${x_1}\bullet ({x_2}*y).$$
\item $K \cap \sett{(k_{x_1},k_{x_2},s,k_y,t)\in K_{x_1}\times K_{x_2}\times [0,1]\times K_y \times [0,1]}{s=1}$ is $$\pow{{x_2}(y-1)}({x_1}*y)\bullet {x_2}$$
\item $K \cap \sett{(k_{x_1},k_{x_2},s,k_y,t)\in K_{x_1}\times K_{x_2}\times [0,1]\times K_y \times [0,1]}{s+t=1}$ is $$({x_1}\bullet {x_2})*y$$
\end{romlist}

then $$
\partial K=K_{\partial y,x_1,x_2}+\pow{y}K_{y\partial {x_1},{x_2}}+\pow{y+{x_1}}K_{y,{x_1},\partial {x_2}}+$$
$$K_{(x_1\bullet
x_2)*y}-K_{x_1\bullet(x_2*y)}-\pow{x_2}K_{({x_1}*y)\bullet
{x_2}}.$$

$\partial \varphi$ is the restriction of
$\varphi$ to $\partial K$. Thus $$(\partial
\varphi -\varphi
\partial)(x_1 \otimes x_2 \otimes y)=(x_1 \bullet x_2)\ast
y -x_1 \bullet (x_2 \ast
y)-\pow{|x_1|(|y|+1)}(x_1 \ast y) \bullet x_2$$
which completes the proof of $(2)$. (See also
\cite{Ger}).
\end{proof}

By \coryc{jac} and \lemc{star} we have,

\start{theo}{gersten} The loop product $\bullet$
with the loop bracket $\{,\}$ makes the loop
homology into a Gerstenhaber algebra, namely:
\begin{numlist}
\item The loop product $\bullet$ defines a graded commutative, associative algebra.
\item $\{,\}$ is a Lie bracket of degree $1$, which means
 that for each $a,b,c \in \HH_\ast$ \begin{romlist}
\item $\{a,b\}=-\pow{(|a|+1)(|b|+1)}\{b,a\}$
\item  $\{a,\{b,c\}\}=\{\{a,b\},c\}+\pow{(|a|+1)(|b|+1)}\{b,\{a,c\}\}$
\end{romlist}
\item $\{a,b\bullet c\}=\{a,b\}\bullet c+\pow{|b|(|a|-1)}b\bullet\{a,c\}$.
\end{numlist}
\end{theo}

\section{The $\Delta$ operator}

Now we consider the degree $+1$ operation on the
chains of the loop space, $$\rightarrow
\flecha{\ls_i}{\lo}{\ls_{i+1}}\rightarrow$$ given
by the circle action on $\Map(\SI, M)$. It can be
defined in the following way: If
$\map{x}{K_x}{\ls}$ is an $i$-chain then
$\map{\lo(x)}{\SI \times K_x}{\ls}$ is the $i+1$
chain such that for each $(s,k_x)\in \SI \times
K_x$, $\lo(x)(s,k_x)(\gamma)=x(k_x)(\gamma+s)$.

Since $\lo$ commutes with the $\partial$ operator
on chains, it passes to the loop homology, the
homology of the free loop space, inducing a
degree $+1$ operator $\Delta$. Moreover, if $x$
is an $i$-chain and $k \ge 1$ then $(\lo)^k(x)$
has always geometric dimension $i+1$. Therefore
we obtain

\start{prop}{ll}
$\map{\Delta}{\HH_\star}{\HH_\star}$ is a degree
$+1$ operator and $\Delta \circ \Delta=0$.
\end{prop}

We want to study how $\Delta$ interacts with the
above structure $\bullet$ and $\{,\}$. In order
to do it, we need to define two auxiliary degree
$+1$ operators on $\ls_\ast$ $$\map{\lo_1,
\lo_{2}}{\ls_\ast}{\ls_\ast}.$$ Let
$\map{x}{K_x}{\ls}$ be a $k$-chain. Then
$$\map{\lo_1(x)}{[0,\frac{1}{2}]\times K_x}{\ls},
\qquad \map{ \lo_{2}(x)}{[\frac{1}{2},1]}{\ls},$$
are the $k+1$-chains defined by\ $$
\lo_1(x)(s,k_x)(\gamma)=x(k_x)(\gamma+s) \mbox{
and }
\lo_{2}(x)(s,k_x)(\gamma)=x(k_x)(\gamma+s).$$
Hence, $\lo=\lo_1+\lo_{2}$

The transversally defined map $$\flecha{\ls_\star
\otimes \ls_\star}{}{\ls_\star},$$
$$\flecha{x\otimes y}{}{x\bullet \lo y}$$ is a
chain map because it is composition of
(transversally defined) chain maps. On the other
hand, using \lemc{lemma2} one can prove that the
map $$\flecha{\ls_\star \otimes
\ls_\star}{}{\ls_\star},$$ $$\flecha{x\otimes
y}{}{\pow{|x|}\lo_{2}(x\bullet y)-x\ast y}$$ is
also a chain map. Moreover, these two chain maps
are chain homotopic, as is shown in the next
lemma.

\start{lem}{lemma4} For $x,y \in \ls_\star$, $$
\pow{|x|}\lo_{2}(x\bullet y)-x\ast y \simeq x
\bullet \lo y $$
\end{lem}
\begin{proof} First, the idea of the proof: consider the chain  operation
$$\map{\varphi}{\ls_\star \otimes
\ls_\star}{\ls_\star}$$ where the loop of $x$ is
attached to any point of the loop of $y$ and one
goes around a part of the loop of $y$, starting
at any point between the marked points of $y$ and
$x$ and ending where $x$ is attached, then goes
around the loop of $x$ and finally, around the
rest of $y$. (See \figc{77}).

\begin{figure}[h]
\begin{center}
\epsfxsize=10cm \epsffile{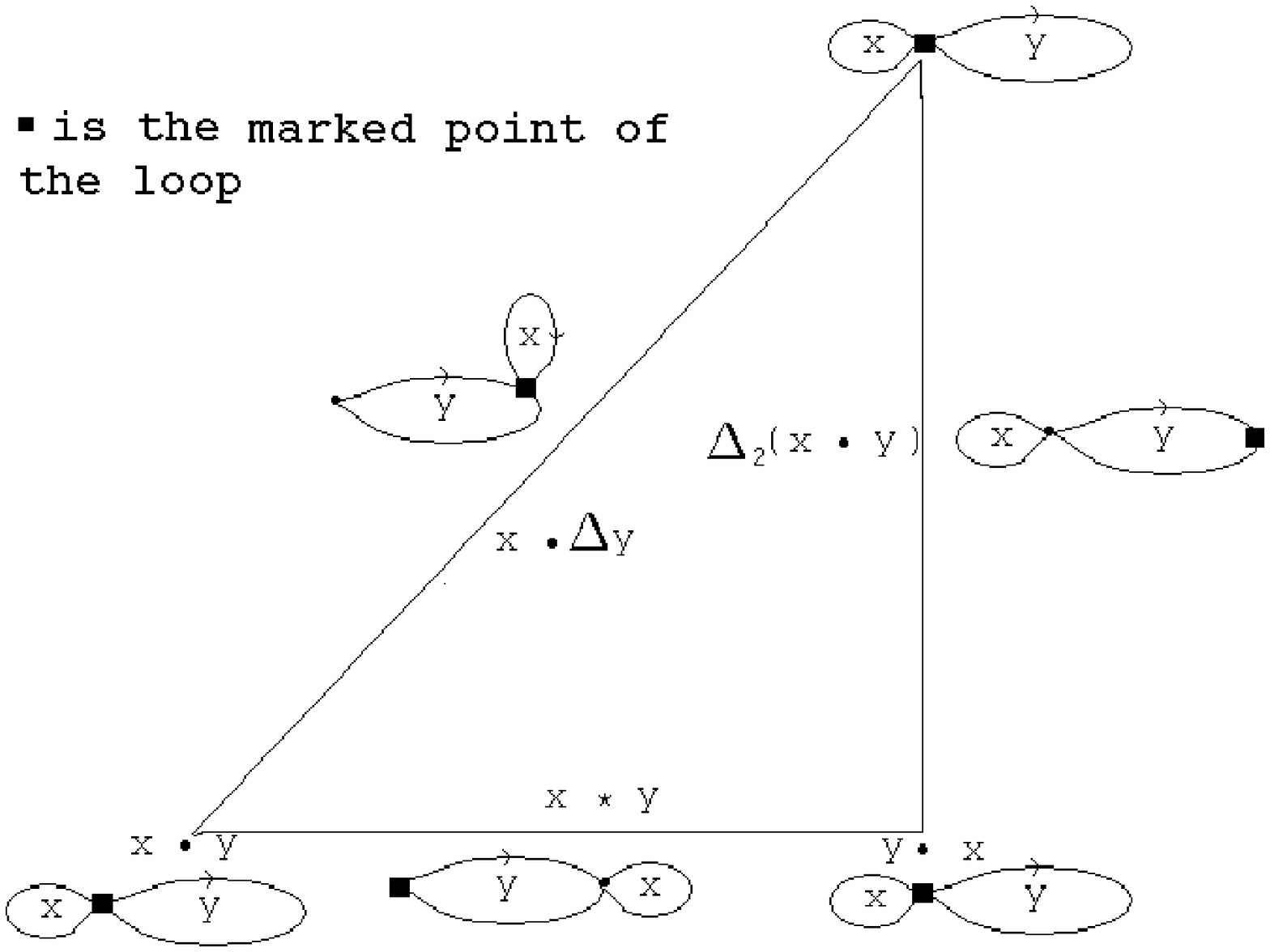}
\end{center}
\caption{Proof of \lemc{lemma4}}\label{77}
\end{figure}

More precisely,

Let $\map{x}{K_x}{\ls}$, $\map{y}{K_y}{\ls}$ be
two cells of $\ls_\star$. Consider $K_{x \ast
y}$, the parameter space of $x \ast y$ and set $$
K=\sett{(t,k_x,s,k_y)\in [0,1]\times K_{x \ast y}
}{ 0 \le t \le s \le 1} $$ and define a chain
$\map{\varphi}{K}{\ls}$ as follows
\begin{eqnarray*}
\varphi(t,k_x,s,k_y)(\gamma) &=
&\left\{\begin{array}{ll} y(k_y)(2\gamma+t)&
\mbox{if $\gamma \in [0,\frac{s-t}{2}]$ ,}\\
x(k_x)(2\gamma-s+t)& \mbox{if $\gamma \in
[\frac{s-t}{2},\frac{s-t+1}{2}]$,}\\
y(k_y)(2\gamma+t)& \mbox{if $\gamma \in
[\frac{s-t+1}{2},1]$.}
\end{array}
\right.
\end{eqnarray*}
Observe that $x$ is attached to $y$ at the image
of $s$ by $y$ and the image of $t$ by $y$ is the
marked point of the resultant loop.

Proceeding in an analogous way as we did in the
proof of \lemc{lemma2}, we obtain $$ (\partial
\varphi-\varphi\partial)(x \otimes
y)=\pow{|x|}\lo_{2}(x\bullet y)-x\ast y - x
\bullet \lo y $$ as desired.
\end{proof}

\start{cory}{delta} The loop bracket $\{,\}$ on
the loop homology is the deviation of $\Delta$
from being a derivation of the loop product. In
other words, for $a, b \in \HH$
$$\{a,b\}=\pow{|a|}\Delta(a \bullet
b)-\pow{|a|}\Delta a \bullet b-a \bullet \Delta b
$$
\end{cory}

\start{theo}{BV} The loop product $\bullet$ and
the operator $\Delta$ make the loop homology into
a Batalin Vilkovisky algebra, namely:
\begin{numlist}
\item $\bullet$ is a graded commutative associative algebra.
\item $\Delta \circ \Delta=0$.
\item $\pow{|a|}\Delta(a \bullet b)-\pow{|a|}\Delta a \bullet b-a \bullet \Delta b$ is a derivation of each variable.
\end{numlist}
\end{theo}
\begin{proof} By \theoc{gersten},\propc{ll} and \coryc{delta}.
\end{proof}

\start{rem}{get} The alternative definition of a
Batalin Vilkovisky algebra as a graded
commutative algebra $(A,\cdot)$ with a degree
$+1$ operator $\map{\Delta}{A}{A}$ such that
$\Delta \circ \Delta=0$ and for each $a,b,c \in
A$ $$ \Delta(a \cdot b \cdot c)= \Delta(a \cdot
b)\cdot c+\pow{|a|}a \cdot \Delta(b \cdot
c)+\pow{(|a|-1)b}b \cdot \Delta(a\cdot c)$$$$
-\Delta(a)\cdot b\cdot c-\pow{|a|}a \cdot
\Delta(b)\cdot c-\pow{|a|+|b|}a\cdot b
\cdot\Delta(c) $$ can be found in Getzler
\cite{Get}.
\end{rem}

\section{The String Bracket}\label{3}

Using the circle action on the free loop space,
we can define the equivariant homology
$\hh_\star$ of the entire loop space. We could
describe this as the ordinary homology of the
quotient space $\s$ of general smooth mappings
$(\SI, M^d)$ by the circle action of rotation in
the domain circle $\SI$. The space $\s$ can be
viewed as the space of all general smooth closed
curves in $M$. Thus we refer to the equivariant
homology of the mapping space $\Map(\SI,M)$ as
{\em string homology}. The circle fibration $$
\flecha{\SI}{}{\mbox{(typical
loops)}}\longrightarrow{\s}\mbox{=string space}
$$ leads to an exact sequence (geometric grading)
$$ \longrightarrow
\flecha{\HH_i}{\erase}{\hh_i}\hflecha{c}{\hh_{i-2}\hflecha{\mark}{\HH_{i-1}}\longrightarrow
\dots} $$ where $\erase$ forgets the marked point
of each member of a family of loops, $\mark$
places a mark on each string in a family in all
possible positions (and $c$ is defined by cap
product with the characteristic class of the
circle bundle above).

The operator $\Delta$ above is the composition
$\mark \circ \erase$. The composition $\erase
\circ \mark$ on homology is zero, as part of the
exactness above.

Any operation $\flecha{\HH_\star^{\otimes
k}}{\wt{\sigma}}{\HH_\star}$ given by composition
of $\bullet$ and $\Delta$ yields an operation
$\flecha{\hh_\star^{\otimes
k}}{\sigma}{\hh_\star}$, $\sigma=\erase \circ
\wt{\sigma} \circ \mark^{\otimes k}$.

In particular, taking $\wt{\sigma} = \bullet$ and
adding a sign, gives the binary operation
$\flecha{\hh_\star \otimes
\hh_\star}{[,]}{\hh_\star}$ called the {\em
string bracket}, $$ [a,b]=\pow{|a|}\erase(\mark
(a) \bullet \mark (b)). $$ where
$|a|=\mbox{dimension } a -d$.

\start{theo}{gol} String homology with the string
bracket, $(\hh_\star,[,])$ is a graded Lie
algebra of degree $(2-d)$ for the geometric
grading.
\end{theo}
\begin{proof} By \theoc{lemma1}, $\bullet$ is graded commutative. So, since $\mark$ has degree $+1$,
$$ [a,b]=\pow{|a|+(|a|+1)(|b|+1)}\erase(\mark (b)
\bullet \mark (a))=-\pow{|a|.|b|}[b,a]. $$ To
prove Jacobi, replace $a$ (resp. $b$, $c$) by
$\mark(a)$ (resp. $\mark(b)$, $\mark(c)$) in the
Leibniz property $(3)$ of \theoc{gersten} and
apply $\erase$ to both sides of the equation to
obtain $$\erase\left(\{\mark(a),\mark(b)\bullet
\mark(c)\}-\{\mark(a),\mark(b)\}\bullet
c-\right.$$$$
\left.\pow{|a|(|b|+1)}\mark(b)\bullet\{\mark(a),\mark(c)\}
\right)=0. $$ Now, use \theoc{BV} and the fact
that $\mark \circ \erase=\Delta$ to replace in
the above equation each of the brackets $\{x,y\}$
by the formula $\powb{c}(\mark \circ
\erase)(c\bullet d)-\powb{c}(\mark \circ
\erase)(c)\bullet d-c\bullet (\mark \circ
\erase)(d).$ Since $\erase \circ \mark=0$ we
cancel the terms where $\erase \circ \mark$
appear and so we obtain $$
\erase\left(-\mark(a)\bullet\mark\erase(\mark(b)\bullet
\mark(c)\right)-\pow{|a|+1}\mark\erase\left(\mark(a)\bullet
\mark(b)\right)\bullet \mark c$$$$
-\pow{|a|(|b|+1)+|a|+1}\mark(b)\bullet\mark
\erase\left(\mark(a)\bullet \mark(c))\right)=0 $$

Now, replacing in the above formula each
occurrence of $\erase(\mark(d),\mark(e))$ by
$\powb{d}[d,e]$ yields $$
-\pow{|a|+|b|}[a,[b,c]]+\pow{|a|+|b|}[[a,b],c]+\pow{|a||b|+|a|+|b|}[b,[a,c]]=0.$$
Hence,
$$[a,[b,c]]=[[a,b],c]+\pow{(|a||b|)}[b,[a,c]].$$
\end{proof}

By the above procedure, considering $\erase \circ
\wt{\sigma} \circ \mark^{\otimes k}$ we define
operations $$ \map{\bar{m}_k}{\hh^{\otimes
k}}{\hh}, \qquad k=2,3,4, \dots $$ by the formula
$$ \flecha{\hh^{\otimes k}}{M^{\otimes
k}}{\HH^{\otimes k}}\hflecha{(\bullet)^{\otimes
(k-1)}}{\HH}\hflecha{\erase}{\hh}. $$ Note that
if we shift the grading on string homology by
$(-d+1)$ from its geometric grading the degree of
each $\bar{m}_k$ becomes
$+1(=-(k-1)d+k-k(-d+1)+(-d+1))$.

Now extend each $\bar{m}_k$ to a coderivation
$m_k$ to $\Lambda \hh_\star$. Here, $\Lambda
\hh_\star$ is the free graded commutative
coalgebra on $\hh_\star$ with the new {\em
algebraic grading} and with the coalgebra
structure which dualizes the usual $\wedge$
algebra structure on $\Lambda \hh^\star$. Thus
$m_k$ is the unique operation whose dual
operation is the (unique) derivation on $\Lambda
\hh^\star$ extending the operator dual to
$\bar{m}_k$.

One knows the Jacobi identity for the string
bracket $\bar{m}_2$ is equivalent to the relation
$m_2 \circ m_2=0$ for the associated coderivation
$m_2$.

The Jacobi relation for $\bar{m}_2$ generalizes
to the entire collection $\{\bar{m}_2, \bar{m}_3,
\dots\}$ in the following way.

\start{theo}{generalized} The associated
coderivations $\{m_2, m_3, \dots\}$ of the free
commutative coalgebra $\Lambda \hh_\star$ on the
string homology $\hh_\star$ satisfy:
\begin{romlist}
\item $m_k \circ m_k=0$, for $k=2,3,4\dots$.
\item $m_k \circ m_r + m_r \circ m_k=0$ for $k,r=2,3,4,\dots$.
\end{romlist}
\end{theo}
\begin{proof}
We only have to show the equations in the case
when the range of the commutator being studied
lies in monomial degree one. Since the commutator
of coderivations is a coderivation, this is
enough to show it is identically zero.

We illustrate the proof of $(ii)$ for $k=3$,
$r=2$. We have four families of closed curves
$A_1, A_2, A_3, A_4$. Then $(m_2 \circ m_3)(A_1
\wedge A_2 \wedge A_3 \wedge A_4)$ can be viewed
a sum of twelve terms, each of them labeled with
$$(\{A_{i_1},A_{i_2},A_{i_3}\}, A_{i_1}),$$ where
$\{A_{i_1},A_{i_2},A_{i_3}\}$ runs over all
possible choices of three families from
$\{A_1,A_2,A_3,A_4\},$ with a preferred element.

Analogously, $(m_2 \circ m_3)(A_1 \wedge A_2
\wedge A_3 \wedge A_4)$ can be viewed as a sum of
twelve terms, each of them labeled with
$$(\{A_{l_1},A_{l_2}\}, A_{l_1}),$$ where
$\{A_{l_1},A_{l_2}\}$ runs over all possible
choices of two families with a preferred element.

A correspondence of the two sets of labels is
given by the map $$(\{A_{i_1},A_{i_2},A_{i_3}\},
A_{i_1})\longrightarrow (\{A_{i_1},A_{n}\},
A_{i_1}) $$ where $A_n$ is the only family not in
$\{A_{i_1},A_{i_2},A_{i_3}\}$.

Now, we will see that corresponding pair of terms
appear in $m_3\circ m_2+m_2\circ m_3$ with
different sign.

Consider, for instance, the corresponding pair of
terms labeled by $$(\{A_{1},A_{2},A_{3}\},
A_{3})\mbox{ and } (\{A_{3},A_{4}\}, A_{3}). $$

For simplicity, let us denote by $A_i$ the
parameter space of the family $A_i$. We can
assume that each $A_i$ is a cell, and that
$\mark(A_i)$ is a map $$\map{\mark(A_i)}{\SI
\times A_i}{\ls}$$ where
$$\mark(A_i)(s,a)(\gamma)= A_i(a)(s+\gamma).$$
The parameter space of the term labeled
$(\{A_{1},A_{2},A_{3}\}, A_{3})$ is $K_\varphi$,
the preimage of $\diag M^3 \times \diag M^2$
under the map $$ \flecha{\SI \times \SI \times
A_1 \times \SI \times A_2 \times \SI \times A_3
\times \SI \times A_4}{\varphi}{M^3 \times M^2},
$$ given by
$\varphi(t,s_1,k_1,s_2,k_2,s_3,k_3,s_4,k_4)=$
$$\left((A_1(k_1)(s_1),A_2(k_2)(s_2),A_3(k_3)(s_3)\right),
\left(A_3(k_3)(t),A_4(k_4)(s_4))\right)$$

The parameter space of the term labeled with
$(\{A_{3},A_{4}\}, A_{3})$, $K_\psi$ is the
preimage of $\diag M^3 \times \diag M^2$ under
the map $$ \flecha{\SI \times \SI \times A_1
\times \SI \times A_2 \times \SI \times A_3
\times \SI \times A_4}{\psi}{M^3 \times M^2}, $$
given by
$\psi(s_1,t,k_1,s_2,k_2,s_3,k_3,s_4,k_4)=$
$$\left((A_1(k_1)(s_1),A_2(k_2)(s_2),A_3(k_3)(s_3)\right),
\left(A_3(k_3)(t),A_4(k_4)(s_4))\right).$$

Over each point of these parameter spaces, the
loops of the terms labeled with
$(\{A_{1},A_{2},A_{3}\}, A_{3})$ and
$(\{A_{3},A_{4}\}, A_{3})$ are as in \figc{term}.

\begin{figure}[h]
\begin{center}
\epsfxsize=8cm \epsffile{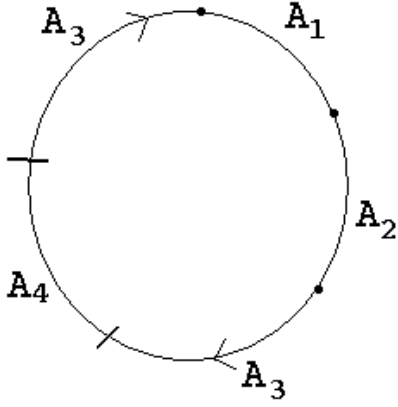}
\end{center}
\caption{The loop of the terms labeled with
$(\{A_{1},A_{2},A_{3}\}, A_{3})$ and
$(\{A_{3},A_{4}\}, A_{3})$}\label{term}
\end{figure}

Observe that the only difference between
$\varphi$ and $\psi$ is that $s_1$ and $t$ are
interchanged. This produces the difference of
sign.

Now, we prove $(i)$. We illustrate for $k=3$.
There will be five families, $A_1, A_2, A_3, A_4,
A_5$. $(m_3 \circ m_3)(A_1\wedge A_2\wedge
A_3\wedge A_4\wedge A_5)$ is the sum of thirty
terms, each of them labeled with
$$(\{A_{i_1},A_{i_2},A_{i_3}\},A_{i_1}) $$ where
$(\{A_{i_1},A_{i_2},A_{i_3}\}$ runs over all
subsets of three elements of $\{A_1, A_2, A_3,
A_4, A_5\}$.

We group these terms in pairs with the following
correspondence
$$(\{A_{i_1},A_{i_2},A_{i_3}\},A_{i_1})\leftrightarrow
(\{A_{i_1},A_{j_1},A_{j_2}\},A_{i_1}) $$ where
$\{A_{j_1},A_{j_2}\}=\{A_1, A_2, A_3, A_4, A_5\}
\setminus \{A_{i_1},A_{i_2},A_{i_3}\}.$ As in the
proof of case $(ii)$, we can see that the
parameter spaces corresponding to pairs cancel,
so $(i)$ holds.
\end{proof}

\start{cory}{uncountable}There exists an
uncountable family $\{\delta_\Lambda\}$ of
$\mathrm{Lie}_\infty$ structures on the string
homology. Namely, for each $\Lambda \subset
\{2,3,\dots\}$, $$\Smap{\delta_\Lambda}{\Lambda
\hh_\star}\mbox{ defined as }
\delta_\Lambda=\sum_{\lambda \in
\Lambda}m_\lambda $$ is a coderivation which
satisfies $\delta_\Lambda \circ
\delta_\Lambda=0$.
\end{cory}
\begin{proof}
By \theoc{generalized} each $\delta_\lambda$ is a
coderivation of $\Lambda \hh_\star$ of square
zero and so determines (by definition) a
$\mathrm{Lie}_\infty$ structure on $\hh_\star$.
\end{proof}

\section{Examples}\label{examples}

\start{ex}{gold} Let us unwrap the string bracket
for $M^d$ when $d=2$. If the genus is greater
than one, the equivariant homology of the loop
space is concentrated in dimension zero (except
for the component of the trivial loop whose
higher homology we ignore at the moment).

The zeroth equivariant homology group, $\hh_0$,
is the vector space with basis $\hat{\pi}$, where
$\hat{\pi}$ denotes the set of free homotopy
classes of loops in $M$.

Let us calculate the string bracket of two
elements $a, b \in \hh_0$. Let $\map{\varrho,
\sigma}{\SI}{M}$ be two loops such that their
free homotopy classes are $a$ and $b$
respectively.

We consider g two elements $x, y$ of $\ls_{-1}$
determined by two maps
$$\flecha{\SI}{x,y}{\Map(\SI, M^2)}$$ given by
$$x(s)(\gamma)=\varrho(s+\gamma) \qquad
y(t)(\gamma)=\sigma(t+\gamma).$$ Hence, $\mark
a=\bar{x}$ and $\mark b=\bar{y}$, where $\bar{x}$
denotes the homology class of a cycle $x$.

Assume that $x$ and $y$ are transversal.

Thus, $K_{x \bullet y}=\sett{(s,t) \in \SI \times
\SI}{\varrho(s)=\sigma(t)}$. By transversality,
the set of intersection points of the curves
$\varrho$ and $\sigma$, $\varrho \cap \sigma$ is
finite. Observe that it is precisely
$\sett{\varrho(s)}{(s,t) \in K_{x \bullet b}}$.

The orientation of $K_{x \bullet y}$ is given by
endowing each of its points with a sign. For each
$(s,t) \in K_{x \bullet y}$, this sign is
precisely, the intersection index
$\epsilon(\varrho,\sigma, \varrho(s))$ of the
loops $\varrho$ and $\sigma$ at $\varrho(s)$.
Then $$\mark a \bullet \mark b=\sum_{p \in
\varrho \cap \sigma}\epsilon(\varrho, \sigma , p)
(\varrho \sharp \sigma)_p,$$ where $(\varrho
\sharp \sigma)_p$ denotes the free homotopy class
which contains the loop product of $\varrho$ and
$\sigma$ at $p$. Hence, $$[a,b]=\erase(\mark a
\bullet \mark b)=\sum_{p \in \varrho \cap
\sigma}\epsilon(\varrho, \sigma , p)
\erase((\varrho \sharp \sigma)_p),$$

It seems remarkable that this formula is well
defined in the vector space of components, it is
skew symmetric, and satisfies Jacobi. Thus the
string bracket becomes for $d=2$ the formula
discovered by Wolpert \cite{Wol} and Goldman
\cite{Gol} which plays a role in the symplectic
structure of the Techmuller space
\cite{Wol}(because of Thurston's earthquakes) and
the symplectic structure on the flat $G$-bundles
for $G$ a compact semisimple Lie group
\cite{Gol}.

Indeed, trying to understand and generalize
Goldman's work \cite{Gol} lead us to the general
theory above.
\end{ex}

\start{ex}{3sphere} Now take $d=3$ and consider a
possibly twisted circle bundle over $M^3$ with
base a surface $F$ of genus greater than one. It
is interesting for the string bracket of $M^3$ to
consider equivariant homology (i.e., string
homology) in dimension one. Cycles are generated
by maps of torii into $M^3$. The projection of
the torus to $F$ is homotopic to a circle. So we
only get interesting examples of torii in $M^3$
by taking all the fibers in $M^3$ over some
circle in $F$. Two of these torii $A$ and $B$ can
be put in transversal position by putting their
projections $a$ and $b$ in $F$ in transversal
position. We see that the string bracket of $A$
and $B$ in $M^3$ is the lift of the string
bracket of $a$ and $b$ in $F$ (as described in
\exc{gold}). So for these $3$-manifolds the
string bracket is just as non trivial as the
string bracket on surfaces.

A similar discussion applies to Seifert
fibrations over surfaces.
\end{ex}

\section{Loop product and cap product}\label{cup}

The cohomology algebra of a space (with cup
product) acts on the homology of a space (called
cap product) via the duality formula $$ <a \cap
x, b>=<a\cup b,x> $$ where $<,>$ is the dual
pairing between homology and cohomology, $\cup$
is cup product, and $\cap$ is cap product. To
relate this structure to our loop product
$\bullet$, consider the diagram $$
\Map(\SI,M)\stackrel{\mbox{c}}{\longleftarrow}\flecha{\Map(\mbox{figure
eight}, M)}{\mbox{i}}{\Map(\SI,M) \times
\Map(\SI,M)} $$ where $c$ denotes the composition
of loops and $i$ is the natural inclusion.
\start{defi}{pair} A pair of cohomology classes
$(a,A)$ (in the appropriate spaces) is called a
{\em compatible pair} if $i^\star A=c^\star a$.
\end{defi}

\start{theo}{compatible} For each $x, y$ homology
classes and compatible pair of classes $(A,a)$ $$
\mbox{loop product }( A \cap(x \otimes y))=a
\cap(x \bullet y). $$
\end{theo}

\begin{proof} The loop product is the composition of intersection with image $i$ (which is represented as a codimension $d$ submanifold being the transverse image of the diagonal under the map
$$\flecha{\Map(\SI,M) \times
\Map(\SI,M)}{\mbox{\scriptsize marked points}}{M
\times M})$$ with the induced transformation of
$c$ in homology. The first process commutes with
cap product with first $A$ then with its
restriction to $\Map(\mbox{figure eight}, M)$.
The second process commutes with capping with a
class in $\Map(\SI, M)$ or with its pull back via
$c$ to $\Map(\mbox{figure eight}, M)$
\end{proof}

\section{Appendix 1: $\Sc^2$ and other simply connected
manifolds}\label{two sphere}

Fibrations such as $$\mbox{based loop maps on
$\Sc^2$}\longrightarrow \mbox{all loops in
$\Sc^2$} \longrightarrow \Sc^2$$ have algebraic
models. For example, $\Sc^2$ and based loops on
$\Sc^2$ are modeled by $\flecha{(0;x;y
\dots)}{d}{(0;0;x^2)}$ and
$\flecha{(x;\bar{y})}{d}{(0;0)}$ respectively.
The notation gives the generators in degrees
$1;2;3;\dots$ respectively and what the
differentials are. The model is the free
commutative algebra with those generators
provided with a derivation $d$ of square zero and
degree $1$ with the specified values.

The total space of the fibration has a model
$\flecha{(\bar{x};\bar{y},x;y)}{{d}}{(0;
-2x\bar{x},0;x^2)}.$

There is also a derivation $\Delta$ of degree
$-1$ and square zero given by
$\flecha{(\bar{x};\bar{y};x;y)}{\Delta}{(0;0;\bar{x};\bar{y})}$
and $d\Delta+\Delta d=0$ is true (and in fact
determines $d$ given $\Delta$ and $d$ on the base
of the fibration.)

The models are cochain models. The obvious maps
serve as cochain maps corresponding to the maps
between space. The cohomology and induced
transformations are derived accordingly
\cite{Sullivan}.

The equivariant cohomology or string cohomology
has model obtained by adding one closed generator
$u$ of degree $2$
$$\flecha{(\bar{x};\bar{y},x,u;y)}{\bar{d}}{(0;
-2x\bar{x},\bar{x}u,0;x^2+\bar{y}u)}$$ where
$\bar{d}u=0$ and $\bar{d}z=dz+ (\Delta z)u$
determines the rule for the other generators $z$.

Calculating with these models we find the Betti
numbers of the loop homology are all $1$ in each
dimension $1,2,3,\dots$ and the image of based
loop homology in free homology is in degrees
$1,3,5,\dots$. Thus exactness for $\Sc^2$ of
$$\mbox{based loop homology}
\stackrel{\mbox{\scriptsize
inclusion}}{\longrightarrow} \flecha{\mbox{loop
homology}}{\cap}{\mbox{based loop homology}}$$
implies $\cap$ has non trivial image in the even
degrees $0,2,4,\dots$. This shows {\em the loop
product is non-trivial for $\Sc^2$}.

Calculating further, one finds $\Delta$ is an
isomorphism on loop homology from odd to even
dimensions and is zero in even dimensions. Since
$\Delta=\mark \circ \erase$, one finds that
$\erase$ is non zero in odd dimensions and
$\mark$ is non zero starting in odd dimensions.
Since $\erase \circ \mark=0$, one finds $\erase$
is zero in even dimensions and $\mark$ is zero
starting in even dimensions. {\em Thus the string
bracket $[x,y]=\pm \erase(\mark(x)\bullet
\mark(y))$ is zero in all dimensions for
$\Sc^2$.}

\subsection{General simply connected manifolds}

This pattern works in general for simply
connected manifold to describe the algorithm for
calculating the rational homology of these
spaces:

If $M$ {\em has minimal model }(where each
differential is quadratic $+$ higher order terms)
$$\flecha{(0;x_1,x_2,\dots;y_1,y_2,\dots)}{d}{(d
x_1,d x_2,\dots;d y_1,d y_2,\dots)}$$ {\em then
the based loops on $M$ has a model }
$$\flecha{(\bar{x_1},\bar{x_2},\dots;\bar{y_1},\bar{y_2},\dots)}{d}{(0,0,\dots0;0,\dots)}$$
and {\em the free loop space with operations
$\Delta$ and $d$ of degree $-1$ and $1$
satisfying $d\Delta+\Delta d=0$ has model with
generators }
$$(\bar{x_1},\bar{x_2},\dots;\bar{y_1},\bar{y_2},\dots,x_1,x_2,\dots;y_1,y_2,\dots)
$$ $\Delta$ is defined by
$\flechi{x_i}{}{\bar{x}_i}$,
$\flechi{y_i}{}{\bar{y}_i}$ and
$\flechi{\bar{x}_i,\bar{y}_i}{}{0}$. Then $d$ is
defined so that $d \Delta+\Delta d=0$ and $d$ is
given as before on the $x_i, y_i \dots$ coming
from $M^d$. Thus $(d\Delta+\Delta d)x_i=0$
implies $d\bar{x}_i=-\Delta(d x_i)$ can be
calculated since $d x_i$ and $\Delta$ are known.

A {\em equivariant model} is obtained by adding
to a free loop space model one more variable $u$
in degree two with $d u=0$ and $\bar{d}z=d
z+\Delta(z)\cdot u$ for the other generators $z$
\cite{Sullivan}.

\start{rem}{it} In \cite{SV} it was shown the
ranks of the loop homology are unbounded for
simply connected manifolds unless the minimal
model has only one or two generators (like
$\Sc^2$ or $\Sc^3$).
\end{rem}

\section{Appendix 2: $M^3$ and $K(\pi,1)$ manifolds}\label{three}

It is known that any closed $3$-manifold is a
connected sum $M_1 \sharp M_2 \sharp M_3
\sharp\dots \sharp M_n$ along $\Sc^2$ where each
of the $M_i$'s is of one of the following types
\begin{romlist}
\item $\pi_1$ is finite so the universal cover is homotopy equivalent to $\Sc^3$.
\item $M_i$ is $\SI \times \Sc^2$.
\item $\pi_1$ is infinite and the universal cover is contractible \cite{JS}
\end{romlist}

The technique of models plus finite group
invariance can be used to treat the examples of
type $(i)$ and $(ii)$ .

If we treat examples of type $(iii)$ it seems
plausible one could develop an algorithm for the
connected sum using our knowledge of $\Sc^2$ and
free product ideas.

We discuss type $(iii)$ under the hypothesis of
Thurston's geometrization picture.
\begin{numlist}
\item If $M_i$ is closed hyperbolic, each centralizer of a non-zero conjugacy class is infinite cyclic. Thus that component of the free loop space is a homotopy circle. By dimension reasons all loop products between these components are zero. The loop product reduces to the classical intersection product.

\item Otherwise $M$ would be a union along torii of Seifert fibrations over surfaces with boundary and finite volume hyperbolic manifolds with neighborhoods of the cusps deleted.

\item If any non trivial Seifert fibrations are present, we have a rich structure of loop product as described in \exc{3sphere}.

\item Finally, $M$ could be a union along torii of hyperbolic pieces and we haven't analyzed these cases.
\end{numlist}
\subsection{General $K(\pi,1)$ manifolds, $d \ge 3$}

The loop space of $M$ is homotopy equivalent to a
union over conjugacy classes $\alpha$ in $\pi$ of
$K(\pi_\alpha,1)$ homotopy types, where
$\pi_\alpha$ is the centralizer of a
representative of $\alpha$. In particular a
component of $\ls(M)$ is homotopy equivalent to a
covering space of $M$. Thus its homological
dimension is at most $d$. Again for hyperbolic
manifolds, the string bracket is zero for
dimension reasons.

\enddocument